\documentclass[12pt]{article}
\usepackage{amsmath,amsfonts,amscd,amssymb}
\newcommand{\IR}{{\mathbb R}} 

\font\sefont=cmr12 scaled \magstep1 
\def\bs{\bigskip}
\def\noi{\noindent}

\def\ms{\medskip}
\def\ve{\varepsilon}

\def\Strich{\par\ms
            \centerline{\hbox to 2cm {\hrulefill}}\par\bs\noi}

\def\de{\partial}

\begin{document}

\thispagestyle{empty}
\def\dg{DGL($\lambda$)\ }
\def\st{_{\textstyle *}}
\def\eps{\ve}
\def\Tr{\mathop{\rm Tr}}
\def\tr{\mathop{\rm tr}}
\def\inv{{}^{-1}}
\def\id{\rm id}
\def\sms{\smallskip}

\centerline{\bf\Large
Bifurcation of straight-line librations}

\vskip.7cm

\centerline{K. J\"anich}
\vskip.5cm
\centerline{\it Faculty of Mathematics, University of Regensburg,
            D-93040 Regensburg, Germany}
\vskip0.8cm

\centerline{\bf Abstract}

\bigskip

\baselineskip 10pt
\noindent
We study a class of 2-dimensional Hamiltonian systems
$H(x,y,p_x,p_y)=\frac12(p_x^2+p_y^2)\\+V(x,y)$ in which the 
plane $x$=$p_x$=0 is invariant under the Hamiltonian flow, so 
that straight-line librations along the $y$ axis exist, and we
also consider perturbations $\delta H=\delta\cdot F(x,y,p_x,p_y)$ 
that preserve these librations. We describe a procedure for the 
analytical calculation of partial derivatives of the Poincar\'e map. 
These partial derivatives can be used to predict the bifurcation 
behavior of the libration, in particular to distinguish between 
transcritical and fork-like bifurcations, as was mathematically 
investigated in \cite{jaen1} and numerically studied in \cite{BT}.
\baselineskip 14pt

\vskip1.cm
\noi
{\sefont 1 Introduction}

\bs\noi
We study 2-dimensional Hamiltonian systems
\begin{equation}
H(x,y,p_x,p_y)=\textstyle\frac12(p_x^2+p_y^2)+V(x,y)
\end{equation}
with a potential satisfying
\begin{equation}
\frac{\de V}{\de x}(0,y)=0\quad\text{for all}\quad y\in\IR\,.
\end{equation}
Then the $(y,p_y)$-plane $x=p_x=0$ is invariant under the Hamiltonian flow
and thus the system will librate on the $y$-axis. Choose one of
the libration families, parametrized by $\ve:=E-E_0\in I$ for some
fixed energy $E_0$ and a suitable open interval $I$.
The family consists of closed orbits
$\gamma_\ve(t)=(0,y(\ve,t),0,p_y(\ve,t))$, with
$p_y(\ve,t)=\dot y(\ve,t))$, which we let start, say, at their maximal
value of $y$, that is at the point
$(0,y_{\rm max}(\ve),0,0)\in\IR^4$. Let $T(\ve)>0$ denote the
period of $\gamma_\ve$.

For each $\ve$ we may use $p_y=0$ as a Poincar\'e surface of
section PSS at the starting point of the orbit, and the
$(x,p_x)$-plane as the PPSS, the projected Poincar\'e surface of
section. As the canonical coordinates $q$ and $p$ in the PPSS, we
may choose $x$ and $p_x$. Then the Poincar\'e map defines a
symplectic family
\begin{equation}
\begin{array}{lll}
Q&=&Q(q,p,\ve)\\
P&=&P(q,p,\ve)\\
\end{array}
\end{equation}
on an open neighborhood of the $\ve$-interval $0\times0\times I$
in the $(q,p,\ve)$-space $\IR^3$. Note that
\begin{equation}
A:=0\times0\times I=\{(0,0,\ve)\mid \ve\in I\}
\end{equation}
itself is a fixed point branch of this family, and we propose to study
the bifurcations that may occur along this branch.

Among the questions we ask about these bifurcations, there is
one concerning the behavior of a bifurcation under a small
deformation of the Hamiltonian. Let $\delta$ denote a small
deformation parameter and let us add a deformation term
$\delta\cdot F(x,y,p_x,p_y)$ to the Hamiltonian:
\begin{equation}
H(x,y,p_x,p_y,\delta)=\textstyle\frac12(p_x^2+p_y^2)+V(x,y)
+\delta F(x,y,p_x,p_y),
\end{equation}
and let $F$ satisfy the following `libration preserving
condition':
\begin{equation}
\frac{\de F}{\de p_x}(0,y,0,p_y)=0
\quad\text{and}\quad
\frac{\de F}{\de x}(0,y,0,p_y)=0.
\end{equation}
Then for fixed $\delta$, the $(y,p_y)$-plane $x=p_x=0$ will
still be invariant under the Hamiltonian flow, so the system will
still librate on the $y$-axis, and we ask how the bifurcation
behavior found at $\delta=0$ will change if we turn on the
parameter $\delta$.

The symplectic family defined by the Poincar\'e map now depends on
two parameters $\ve$ and $\delta$:
\begin{equation}
\begin{array}{lll}
Q&=&Q(q,p,\ve,\delta)\\
P&=&P(q,p,\ve,\delta).\\
\end{array}
\end{equation}
Think of an arbitrary $\ve_0$ being chosen. We will ask if
the fixed point $(0,0)$, of the undeformed system, is singular
at $\ve=\ve_0$ and if so, what are the properties of the bifurcation
and their behavior under deformations $\delta\neq0$. Using
\cite{jaen1}, the answers to these questions could be read from those 
38 partial derivatives up to third order of $P$ and $Q$
that involve the parameters $\ve$ and $\delta$ at most in first order, at
the single point $(0,0,\ve_0,0)$ --- if we only knew them. The
purpose of the present note is to describe a procedure for the
calculation of these partial derivatives of the Poincar\'e map
from the potential $V(x,y)$ and the deformation term
$F(x,y,p_x,p_y)$.

\vskip1.5cm

\noi
{\sefont 2 Numerical prerequisites}

\bs\noi
To start the procedure, for a given $\ve_0$, we will need to
know first of all the closed orbit $\gamma_{\ve_0}$ itself,
$\gamma_{\ve_0}(t)=(0,y(t,\ve_0),0,p_y(t,\ve_0))$, that is we
have to solve the equation
\begin{equation}
\ddot y+\frac{\de V}{\de y}(0,y)=0
\end{equation}
to the initial condition $y(0)=y_{\rm max}(\ve_0)$ and
$\dot y(0)=0$. The value $y_{\rm max}(\ve_0)$ satisfies
\begin{equation}
V(0,y_{\rm max}(\ve_0))=E_0+\ve_0.
\end{equation}
Within the chosen domain of libration it will be the larger of
the two solutions of this equation and can be determined that
way. The function $y(t):=y(t,\ve_0)$ will be periodic of a period
$T(\ve_0)>0$.

Furthermore, we will have to solve the linearized Hamiltonian
equation along this closed orbit, that is we have to know the
fundamental system $(\xi_1, \xi_2)$ of the linear equation
\begin{equation}
\ddot \xi+\frac{\de^2 V}{\de x^2}(0,y(t))\xi=0
\end{equation}
to the initial condition
\begin{equation}
\left(
\begin{array}{ll}
\xi_1(0)&\xi_2(0)\\
\noalign{\sms}
\dot\xi_1(0)&\dot\xi_2(0)\\
\end{array}
\right)
=
\left(
\begin{array}{lc}
1&0\\
\noalign{\sms}
0&1\\
\end{array}
\right)
\end{equation}
as well as the fundamental system $(\eta_1,\eta_2)$ of
\begin{equation}
\ddot \eta+\frac{\de^2 V}{\de y^2}(0,y(t))\eta=0
\end{equation}
to the initial condition
\begin{equation}
\left(
\begin{array}{ll}
\eta_1(0)&\eta_2(0)\\
\noalign{\sms}
\dot\eta_1(0)&\dot\eta_2(0)\\
\end{array}
\right)
=
\left(
\begin{array}{lc}
1&0\\
\noalign{\sms}
0&1\\
\end{array}
\right).
\end{equation}
To say that these five functions $y(t)$, $\xi_1(t)$, $\xi_2(t)$,
$\eta_1(t)$, $\eta_2(t)$ and their first derivatives
must be `known' means that they
are known numerically on the whole period interval
$[0,T(\ve_0)]$. A computer program implementing the procedure
for the calculation of the partial derivatives of the Poincar\'e
map at $(0,0,\ve_0,0)$ will have to treat them as known
functions. But beyond that no further differential equations
will have to be solved.

\vskip1.5cm

\noi
{\sefont 3 The Poincar\'e map}

\bs\noi
The Poincar\'e map is produced by the Hamiltonian flow. To
facilitate the handling of higher partial derivatives, we will
gradually shift from natural `speaking' notations like
$(x,y,p_x,p_y)$ to a simple enumeration of variables by upper
indices. Partial derivatives by these variables will then be
denoted by corresponding lower indices. We begin by writing
\begin{equation}
\begin{array}{l}
a^1:=x\\
a^2:=y\\
a^3:=p_x\\
a^4:=p_y\\
\end{array}
\end{equation}
for the independent variables in $\IR^4$ and
\begin{equation}
x^i=x^i(t,a^1,a^2,a^3,a^4,\delta)
\end{equation}
with $i=1,\dots,4$, for the components of the Hamiltonian flow at a fixed
$\delta$, with initial conditions $a^1,\dots, a^4$ :
\begin{equation}
x^i(0,a^1,a^2,a^3,a^4,\delta)=a^i.
\end{equation}
For fixed $\ve$ and $\delta$, the starting point in the PSS
corresponding to a given point $(q,p)$ in the PPSS is described
by
\begin{equation}
\begin{array}{l}
a^1=q\\
a^2=y(q,p,\ve,\delta)\\
a^3=p\\
a^4=0,\\
\end{array}
\end{equation}
where the $y$-component is defined implicitly by
\begin{equation}
\frac12p^2+V(q,y(q,p,\ve,\delta))+
\delta F(q,y(q,p,\ve,\delta),p,0)=\ve+E_0
\end{equation}
and $y(0,0,\ve,0)=y_{\rm max}(\ve)$. From this starting point,
the flow line will travel for a time $T=T(q,p,\ve,\delta)>0$
until it hits the PSS $p_y=0$ again, so implicitly this time is
given by
\begin{equation}
x^4(T(q,p,\ve,\delta),q,y(q,p,\ve,\delta),p,0,\delta)=0,
\end{equation}
in the notation (15) of the flow, and by the reference condition
$T(0,0,\ve,0)=T(\ve)$, the period of $\gamma_\eps$. The Poincar\'e
map can now be described as
\begin{equation}
\begin{array}{lll}
Q(q,p,\ve,\delta)&=&x^1(T(q,p,\ve,\delta),q,y(q,p,\ve,\delta),p,0,\delta)\\
P(q,p,\ve,\delta)&=&x^3(T(q,p,\ve,\delta),q,y(q,p,\ve,\delta),p,0,\delta).\\
\end{array}
\end{equation}
Taking partial derivatives, we obtain the the partial
derivatives of the Poincar\'e map in terms of partial derivatives
of the Hamiltonian flow and of partial derivatives of the
flow time function $T(q,p,\ve,\delta)$ and the starting point
function $y(q,p,\ve,\delta)$. That's what we do next.

\vskip1.5cm

\noi
{\sefont 4 Taking derivatives of the Poincar\'e map}

\bs\noi
We are now unifying the notation of the independent variables of
the flow, currently written as $(t,a^1,a^2,a^3,a^4,\delta)$, to
$(a^0,\dots,a^5)$. Derivatives are denoted by lower indices, so
for instance $x^1_{02}$ would mean
\begin{equation}
x^1_{02}=\frac{\de^2 x^1}{\de t\de(a^2)}=
\frac{\de \dot x^1}{\de(a^2)},
\end{equation}
and so on. Next, we write the two components $Q$ and $P$ of the
Poincar\'e map as compositions
\begin{equation}
\begin{array}{l}
Q=x^1\circ Z\\
P=x^3\circ Z,\\
\end{array}
\end{equation}
or
$Q(q,p,\ve,\delta)=
x^1(Z^0(q,p,\ve,\delta),\dots,Z^5(q,p,\ve,\delta))$ and
analogously for $P$,
where the six components $Z^0,\dots,Z^5$ are given, according to
(20), by
\begin{equation}
\begin{array}{l}
Z^0(q,p,\ve,\delta):=T(q,p,\ve,\delta)\\
Z^1(q,p,\ve,\delta):=q\\
Z^2(q,p,\ve,\delta):=y(q,p,\ve,\delta)\\
Z^3(q,p,\ve,\delta):=p\\
Z^4(q,p,\ve,\delta):=0\\
Z^5(q,p,\ve,\delta):=\delta.\\
\end{array}
\end{equation}
Finally, let us also enumerate the
independent coordinates $(q,p,\ve,\delta)$ in the product of the PPSS with the
parameter plane by writing
\begin{equation}
\begin{array}{l}
u^1:=q\\
u^2:=p\\
u^3:=\ve\\
u^4:=\delta,\\
\end{array}
\end{equation}
again denoting partial derivatives by lower indices, like in
$Z^2_3=\frac{\de Z^2}{\de (u^3)}=\frac{\de y}{\de\ve}$. I shall
try to be consistent in using greek letters $\lambda$, $\mu$,
\ldots for the $a$-related indices that run from 0 to 5 and use
roman letters $\ell, m,\ldots=1,2,3,4$ for indices refering to
the $u$-variables.

The partial derivatives we are after, like
$P_{qq\ve}=(x^3\circ Z)_{113}$, can now neatly be written in
terms of partial derivatives of flow, flow time and starting
point functions as
\begin{equation}
\hskip-6pt\begin{array}{lll}
(x^i\circ Z)_\ell&=&x^i_\lambda Z^\lambda_\ell\\
\noalign{\ms}
(x^i\circ Z)_{\ell m}&=&
x^i_{\lambda\mu} Z^\lambda_\ell Z^\mu_m+
x^i_\lambda Z^\lambda_{\ell m}\\
\noalign{\ms}
(x^i\circ Z)_{\ell mn}&=&
x^i_{\lambda\mu\nu} Z^\lambda_\ell Z^\mu_mZ^\nu_n
+x^i_{\lambda\mu}(
Z^\lambda_{\ell}Z^\mu_{mn}+
Z^\lambda_{m}Z^\mu_{n\ell}+
Z^\lambda_{n}Z^\mu_{\ell m})+
x^i_{\lambda} Z^\lambda_{\ell mn}.\\
\end{array}
\end{equation}
That's just the chain rule, so it holds everywhere. But
for our application, we need to know the left hand sides only at
the point $(u^1,u^2,u^3,u^4)=(0,0,\ve_0,0)$, and so on the right
hand side we want to know the $Z^\lambda\st$ at $(0,0,\ve_0,0)$
and the $x^i\st$ at $(T(\ve_0),0,y_{\rm max}(\ve_0),0,0,0)$.

\vskip1.5cm

\noi
{\sefont 5 Derivatives of the starting point function}

\bs\noi
First we will take care of the derivatives
$Z^\lambda\st$. For $\lambda=1,3,4,5$ this is easy,
because for those $\lambda$ the $Z^\lambda\st$ are given by
\begin{equation}
\begin{array}{l}
Z^1_1=Z^3_2=Z^5_4=1\quad\text{and}\\
\noalign{\ms}
Z^1\st=Z^3\st
=Z^4\st=Z^5\st=0\\
\end{array}
\end{equation}
for all other lower indices $*$, in particular for all derivatives
of order $\ge2$. It remains to determine the derivatives
$Z^0\st$ and $Z^2\st$
of $Z^0=T(u^1,\dots u^4)$ and $Z^2=y(u^1,\dots u^4)$.
In the present section we will calculate $Z^2\st$
up to second order.

The defining equation (18) for the
starting point function $y(u^1,\dots,u^4)$
in the $u$-notation becomes
\begin{equation}
\frac12u^2\cdot u^2+V(u^1,y(u^1,\dots,u^4))+
u^4 F(u^1,y(u^1,\dots,u^4),u^2,0)=u^3+E_0.
\end{equation}
As a first step, we will determine the eight derivatives
\begin{equation}
y_1, y_2, y_3\quad\text{and}\quad
y_{11}, y_{12}, y_{22},y_{13},y_{23}
\end{equation}
at $(0,0,\ve_0,0)$. Since here the variable $u^4=\delta$ is not
involved, we may put $u^4=0$.
Taking first derivatives by $u^1$, $u^2$, $u^3$ gives
\begin{equation}
\begin{array}{rll}
V_1+V_2y_1&=&0\\
\noalign{\ms}
u^2+V_2y_2&=&0\\
\noalign{\ms}
V_2y_3&=&1\\
\end{array}
\end{equation}
at all $(u^1,u^2,u^3,0)$. Differentiating further, we get
\begin{equation}
\begin{array}{rll}
V_{11}+V_{12}y_1+(V_{21}+V_{22}y_1)y_1+V_2y_{11}&=&0\\
\noalign{\ms}
V_{12}y_1+V_{22}y_2y_1+V_2y_{12}&=&0\\
\noalign{\ms}
1+V_{22}y_2y_2+V_2y_{22}&=&0\\
\noalign{\ms}
(V_{21}+V_{22}y_1)y_3+V_2y_{31}&=&0\\
\noalign{\ms}
V_{22}y_2y_3+V_2y_{32}&=&0,\\
\end{array}
\end{equation}
again at all $(u^1,u^2,u^3,0)$. Specializing now to
$(0,0,\ve_0,0)$ and using our assumption
$\frac{\de V}{\de x}(0,y)\equiv0$ on the potential, which
implies
\begin{equation}
V_1(0,y)=V_{12}(0,y)\equiv0
\end{equation}
for all $y$, we obtain
\begin{equation}
\begin{array}{rlc}
y_1=y_2=y_{12}=y_{13}=y_{23}&=&\phantom{-}0\\
\noalign{\ms}
\text{and}\quad y_3&=&\displaystyle\phantom{-}\frac1{V_2}\\
\noalign{\ms}
y_{11}&=&\displaystyle-\frac{V_{11}}{V_2}\\
\noalign{\ms}
y_{22}&=&\displaystyle-\frac1{V_2}\\
\end{array}
\end{equation}
at $(0,0,\ve_0,0)$ as some first `end results' on the
derivatives
of the starting point function in terms of
the derivatives $V_2(0,y_{\rm max}(\ve_0))$ and
$V_{11}(0,y_{\rm max}(\ve_0))$ of the potential.

Let us now consider the derivatives $y_4$, $y_{14}$, and
$y_{24}$ involving the deformation parameter $u^4=\delta$.
Differentiating (27) first by $u^4$ and then in addition by
$u^1$ resp. $u^2$ we obtain
\begin{equation}
\begin{array}{rll}
V_2y_4+F&=&0\\
\noalign{\ms}
(V_{21}+V_{22}y_1)y_4+V_2y_{41}+F_1+F_2y_1&=&0\\
\noalign{\ms}
V_{22}y_2y_4+V_2y_{42}+F_2y_2&=&0\\
\end{array}
\end{equation}
at all $(u^1,u^2,u^3,u^4)$. Specializing again to
$(0,0,\ve_0,0)$ and using (31), (32) and the assumption (6) about
$F(x,y,p_x,p_y)$,  we get
\begin{equation}
\begin{array}{rll}
y_{14}=y_{24}&=&\phantom{-}0\\
\noalign{\ms}
\text{and}\quad y_4&=&
\displaystyle-\frac{F(0,y_{\rm max}(\ve_0),0,0)}{V_2(0,y_{\rm max})}.\\
\end{array}
\end{equation}
Note that with (26), (32) and (34) we have determined all
partial derivatives $Z^\lambda_\ell$ and $Z^\lambda_{\ell m}$ at
$(0,0,\ve_0,0)$ for $\lambda\ge1$ and at most one of the indices
$\ell$ and $m$ being 3 or 4. For convenience, let us collect these
results. First for the `easy' $\lambda$'s. Here the first
derivatives are, everywhere:
\begin{equation}
\begin{array}{r|cccc}
Z^\lambda_\ell&{\scriptstyle\ell=1}&{\scriptstyle\ell=2}
&{\scriptstyle\ell=3}&{\scriptstyle\ell=4}\\
\hline
{\scriptstyle\lambda=1}&1&{0}&{0}&{0}\\
{\scriptstyle\lambda=3}&{0}&1&{0}&{0}\\
{\scriptstyle\lambda=4}&{0}&{0}&0&{0}\\
{\scriptstyle\lambda=5}&{0}&{0}&{0}&1\\
\end{array}
\end{equation}
The higher derivatives therefore are zero, in particular
$Z^\lambda_{\ell m}=0$ for $\lambda=1,3,4,5$. For $\lambda=2$ we
have found at $(u^1,u^2,u^3,u^4)=(0,0,\ve_0,0)$:
\begin{equation}
\begin{array}{llc}
Z^2_1&=&\phantom{-}0\\
\noalign{\ms}
Z^2_2&=&\phantom{-}0\\
\noalign{\ms}
Z^2_3&=&\displaystyle\phantom{-}\frac1{V_2}\\
\noalign{\ms}
Z^2_4&=&\displaystyle-\frac{F}{V_2}\\
\noalign{\ms}
Z^2_{11}&=&\displaystyle-\frac{V_{11}}{V_2}\\
\noalign{\ms}
Z^2_{12}&=&\phantom{-}0\\
\noalign{\ms}
Z^2_{22}&=&\displaystyle-\frac{1}{V_2}\\
\noalign{\ms}
Z^2_{13}&=&\phantom{-}0\\
\noalign{\ms}
Z^2_{23}&=&\phantom{-}0\\
\noalign{\ms}
Z^2_{14}&=&\phantom{-}0\\
\noalign{\ms}
Z^2_{24}&=&\phantom{-,}0,\\
\end{array}
\end{equation}
where the derivatives of $V$ have to be taken at
$(x,y)=(0,y_{\rm max}(\ve_0))$, the value of $F$ at
$(x,y,p_x,p_y)=(0,y_{\rm max}(\ve_0),0,0)$.
In the next section, we turn to the
remaining case $\lambda=0$.

\vskip1.5cm

\noi
{\sefont 6 Derivatives of the flow time function}

\bs\noi
Although we are interested in the partial derivatives of
$x^i\circ Z$ for $i=1,3$ only, we will also become involved with the
$p_y$-component $x^4(a^0,\dots,a^5)$ of the flow, because the
defining condition (19) of the flow time function $Z^0$ is
\begin{equation}
x^4\circ Z\equiv0.
\end{equation}
This is also the reason, by the way, why the linear equation
(12) will come up in the calculations. --- From (37), using the
first equation of (25) for $i=4$, we obtain
\begin{equation}
x^4_0Z^0_\ell=-x^4_1Z^1_\ell-x^4_2Z^2_\ell-x^4_3Z^3_\ell
              -x^4_4Z^4_\ell-x^4_5Z^5_\ell.
\end{equation}
The $Z^\lambda_\ell$ on the right hand side are known at
$(0,0,\ve_0,0)$. Correspondingly, the $x^4_\lambda$, for
$\lambda=0,\dots,5$ are meant to be taken at
$(a^0,\dots,a^5)=(T(\ve_0),0,y_{\rm max}(\ve_0),0,0,0)$. We
still have to determine them, but we certainly know $x^4_0$
there, since by the Hamiltonian equations
\begin{equation}
x^4_0(t,0,y_{\rm max}(\ve_0),0,0,0)=
\dot p_y(t)=-\frac{\de V}{\de y}(0,y(t,\ve_0))
\end{equation}
along the closed orbit $\gamma_{\ve_0}$, and therefore
\begin{equation}
x^4_0(T(\ve_0),0,y_{\rm max}(\ve_0),0,0,0)
=-V_2(0,y_{\rm max}(\ve_0)).
\end{equation}
So from (35), (36) and (38) we get
\begin{equation}
\begin{array}{lll}
Z^0_1&=&\displaystyle\phantom{-}\frac1{V_2} x^4_1\\
\noalign{\ms}
Z^0_2&=&\displaystyle\phantom{-}\frac1{V_2} x^4_3\\
\noalign{\ms}
Z^0_3&=&\displaystyle\phantom{-}\frac1{V_2V_2} x^4_2\\
\noalign{\ms}
Z^0_4&=&\displaystyle-\frac{F}{V_2V_2}
x^4_2+\frac1{V_2} x^4_5.\\
\end{array}
\end{equation}
Similarly, from the second equation of (25) for $i=4$ we now
have
\begin{equation}
\begin{array}{lll}
Z^0_{\ell m}&=&\displaystyle\frac1{V_2}\big(
x^4_1Z^1_{\ell m}+
x^4_2Z^2_{\ell m}+
x^4_3Z^3_{\ell m}+
x^4_4Z^4_{\ell m}+
x^4_5Z^5_{\ell m}+
x^4_{\lambda\mu} Z^\lambda_\ell Z^\mu_m
\big)\\
\noalign{\ms}
&=&\displaystyle\frac1{V_2}\big(
x^4_2Z^2_{\ell m}+
x^4_{\lambda\mu} Z^\lambda_\ell Z^\mu_m
\big)\\
\end{array}
\end{equation}
at $(0,0,\ve_0,0)$, with the $x^4\st$ to be taken at
$(T(\ve_0),0,y_{\rm max}(\ve_0),0,0,0)$, as before.
And here we leave it for now:
the $Z^\lambda\st$ on the right
hand side, at $(0,0,\ve_0,0)$, are all known from
(35), (36) and (41). A computer will understand (42) as given, and
for us, there is no point in writing out the formula in great
length before we know more about the $x^4_{\lambda\mu}$, in
particular before we know which of them will vanish anyway.

The
same reasoning applies to the last $Z^\lambda\st$
that are still missing, namely the third derivatives
$Z^2_{\ell m n}$ and $Z^0_{\ell m n}$ of the starting point and
flow time functions, which might be needed in
the computation of the third derivatives of the Poincar\'e map.
In fact they will {\it not} be needed, because they enter the
third equation of (25) for $i=1,3$ with coefficients $x^i_2$ and
$x^i_0$, which will soon be seen to vanish for $i=1,3$. This is just
one of the details of the problem to calculate all the
$x^i\st(T(\ve_0),0,y_{\rm max}(\ve_0),0,0,0)$ that
we need. To this problem our bifurcation analysis is now reduced
and it will be solved in the remaining sections.

\vskip1.5cm

\noi
{\sefont 7 Equations for the flow derivatives}

\bs\noi
The four components of the Hamiltonian flow, see (15) and (16),
are written currently as
\begin{equation}
x^i=x^i(a^0,\dots,a^5)\quad\text{for}\quad i=1,2,3,4
\end{equation}
with $a^0$ denoting the time $t$ and $a^5=\delta$, while
$(a^1,\dots,a^4)$ is the initial point. As a bookkeeping device
we now introduce a `fifth flow component' by
\begin{equation}
x^5(a^0,\dots,a^5):=a^5.
\end{equation}
But we also use $x^1,\dots,x^5$ as the names of the independent
variables of the Hamiltonian vector field
$\vec v=(v^1,\dots,v^4)$, which is then given by
\begin{equation}
\begin{array}{lll}
v^1(x^1,\dots,x^5)&=&\phantom{-}x^3+x^5{} F_3(x^1,\dots,x^4)\\
\noalign{\ms}
v^2(x^1,\dots,x^5)&=&\phantom{-}x^4+x^5{} F_4(x^1,\dots,x^4)\\
\noalign{\ms}
v^3(x^1,\dots,x^5)&=&-V_1(x^1,x^2)-x^5{} F_1(x^1,\dots,x^4)\\
\noalign{\ms}
v^4(x^1,\dots,x^5)&=&-V_2(x^1,x^2)-x^5{} F_2(x^1,\dots,x^4).\\
\noalign{\ms}
\end{array}
\end{equation}
The Hamilton equations become
$\dot x^r(\vec a)=v^r(\vec x(\vec a))$ for
r=1,2,3,4,
and as in (25) we obtain
\begin{equation}
\hskip-6pt\begin{array}{lll}
\dot x^r_\lambda&=&v^r_ix^i_\lambda\\
\noalign{\ms}
\dot x^r_{\lambda\mu}&=&
v^r_{ij}x^i_\lambda x^j_\mu+v^r_ix^i_{\lambda\mu}\\
\noalign{\ms}
\dot x^r_{\lambda\mu\nu}&=&
v^r_{ijk}x^i_\lambda x^j_\mu x^k_\nu
+v^r_{ij}(x^i_{\lambda}x^j_{\mu\nu}+
x^i_{\mu}x^j_{\nu\lambda}+x^i_{\nu}x^j_{\lambda\mu})
+v^r_{i}x^i_{\lambda\mu\nu}\\
\end{array}
\end{equation}
The summation indices $i,j,k$ run from 1 to 5, while any
$\lambda,\mu,\nu\in\{0,1,2,3,4,5\}$ are admitted.
The equations hold everywhere, that is the flow components $x^i$
and their derivatives may be taken at any $\vec a:=(a^0,\dots,a^5)$, the
vector field components $v^r$ and their derivatives then at the
corresponding $\vec x(\vec a))=(x^1(\vec a),\dots,x^5(\vec a))$.
Also note that the time derivative of the flow is the partial
derivative by $a^0$, so on the left hand sides we might have
written $x^r_{0\lambda}$, $x^r_{0\lambda\mu}$,
$x^r_{0\lambda\mu\nu}$ instead of
$\dot x^r_{\lambda}$, $\dot x^r_{\lambda\mu}$,
$\dot x^r_{\lambda\mu\nu}$.

\sms
To apply (25), we will only need to know the
$x^i\st(T(\ve_0),0,y_{\rm max}(\ve_0),0,0,0)$. But in
order to determine these numbers, we will also have to consider
the functions
$x^i\st(t,0,y_{\rm max}(\ve_0),0,0,0)$ on the interval
$[0,T(\ve_0)]$, for which we now introduce the notation
\begin{equation}
x^i\st(t):=
x^i\st(t,0,y_{\rm max}(\ve_0),0,0,0).
\end{equation}
Correspondingly, we write
\begin{equation}
v^r\st(t):=
v^r\st(\vec x(t,0,y_{\rm max}(\ve_0),0,0,0)).
\end{equation}
Then from (47) we get
\begin{equation}
\hskip-6pt\begin{array}{lll}
\dot x^r_\lambda(t)&=&v^r_i(t)x^i_\lambda(t)\\
\noalign{\ms}
\dot x^r_{\lambda\mu}(t)&=&
v^r_{ij}(t)x^i_\lambda(t) x^j_\mu(t)+v^r_i(t)x^i_{\lambda\mu}(t)\\
\noalign{\ms}
\dot x^r_{\lambda\mu\nu}(t)&=&
v^r_{ijk}(t)x^i_\lambda(t) x^j_\mu(t) x^k_\nu(t)\\
\noalign{\ms}
&&+v^r_{ij}(t)\big(x^i_{\lambda}(t)x^j_{\mu\nu}(t)+
x^i_{\mu}(t)x^j_{\nu\lambda}(t)
+x^i_{\nu}(t)x^j_{\lambda\mu}(t)\big)\\
\noalign{\ms}
&&+v^r_{i}(t)x^i_{\lambda\mu\nu}(t)\\
\end{array}
\end{equation}
from which we will now proceed to determine the functions (47).
We only consider $r=1,\dots,4$ as there is no need to write
equations for $x^5\st(t)$, since of course
\begin{equation}
x^5_5\equiv1\quad\text{and}\quad x^5\st(t)\equiv0
\end{equation}
for all other partial derivatives of $x^5(\vec a)=a^5$.

\vskip1.5cm

\newpage

\noi
{\sefont 8 Calculation of the first order flow derivatives}

\bs\noi
Note that all the $v^r\st(t)$ on the right hand sides of (49)
are {\it known} functions in the sense agreed upon in section~2
on the numerical prerequisites, since
\begin{equation}
\vec x(t,0,y_{\rm max}(\ve_0),0,0,0)=
(0,y(t),0,\dot y(t),0),
\end{equation}
which is to be used in (48). In particular, let us tabulate the
first derivatives $v^r_i(t)$:
\begin{equation}
\begin{array}{r|ccccc}
v^r_i(t)&{\scriptstyle i=1}&{\scriptstyle i=2}
&{\scriptstyle i=3}&{\scriptstyle i=4}&{\scriptstyle i=5}\\
\hline
{\scriptstyle r=1}&0&{0}&{1}&{0}&0\\ 
{\scriptstyle r=2}&{0}&0&{0}&{1}&\phantom{-}F_4(t)\\
{\scriptstyle r=3}&{-V_{11}(t)}&{0}&0&{0}&0\\ 
{\scriptstyle r=4}&{0}&{-V_{22}(t)}&{0}&0&-F_2(t)\\
\end{array}
\end{equation}
Here we write $F_4(t):=F_4(0,y(t),0,\dot y(t))$ and
$V_{11}(t):=V_{11}(0,y(t))$ and so on, in line with the notation
introduced in (47) and (48). Note that $F_1(t)=F_3(t)=0$ by our
assumption (6) on $F(x,y,p_x,p_y)$.
As a first consequence of (49), we
see that the $x^r_\lambda(t)$ for $r,\lambda\in\{1,2,3,4\}$ satisfy the
homogeneous linear differential equations
\begin{equation}
\begin{array}{llrll}
\dot x^1_\lambda&-&x^3_\lambda&=&0\\
\noalign{\ms}
\dot x^3_\lambda&+&V_{11}(t)x^1_\lambda&=&0\\
\end{array}
\end{equation}
and
\begin{equation}
\begin{array}{llrll}
\dot x^2_\lambda&-&x^4_\lambda&=&0\\
\noalign{\ms}
\dot x^4_\lambda&+&V_{22}(t)x^2_\lambda&=&0\\
\end{array}
\end{equation}
These are just the first order systems corresponding to (10) and
(12), and the
{\it initial conditions}, as we see from (16), are
\begin{equation}
\begin{array}{lll}
x^r_\lambda(0)&=&1\quad\text{if}\quad r=\lambda,\quad\text{and}\\
\noalign{\sms}
x^r_\lambda(0)&=&0\quad\text{if}\quad r\ne \lambda.\\
\end{array}
\end{equation}
But this shows that the $x^r_\lambda(t)$ for $r,\lambda=1,\dots,4$ are known
functions from the numerical prerequisites, more precisely
\begin{equation}
\left(
\begin{array}{ll}
x^1_1(t)&x^1_3(t)\\
\noalign{\sms}
x^3_1(t)&x^3_3(t)\\
\end{array}
\right)
=
\left(
\begin{array}{ll}
\xi_1(t)&\xi_2(t)\\
\noalign{\sms}
\dot\xi_1(t)&\dot\xi_2(t)\\
\end{array}
\right)\phantom{,}
\end{equation}
and
\begin{equation}
\left(
\begin{array}{ll}
x^2_2(t)&x^2_4(t)\\
\noalign{\sms}
x^4_2(t)&x^4_4(t)\\
\end{array}
\right)
=
\left(
\begin{array}{ll}
\eta_1(t)&\eta_2(t)\\
\noalign{\sms}
\dot\eta_1(t)&\dot\eta_2(t)\\
\end{array}
\right),
\end{equation}
while
\begin{equation}
x^r_\lambda(t)=0\quad\text{for}\quad r+\lambda\quad\text{odd},
\quad\lambda\in\{1,2,3,4\}.
\end{equation}
What about $\lambda=0$ and
$\lambda=5$? We know $x^r_0(t)=\dot x^r(t)=v^r(t)$, so from (45)
we have
\begin{equation}
\begin{array}{lcc}
x^1_0(t)&=&0\\
\noalign{\sms}
x^2_0(t)&=&\dot y(t)\\
\noalign{\sms}
x^3_0(t)&=&0\\
\noalign{\sms}
x^4_0(t)&=&-V_2(t).\\
\end{array}
\end{equation}
Looking at $\lambda=5$, we see from (47) and (52) that (53) is
satisfied also in this case, but (54) has to be replaced by the
{\it inhomogeneous} system
\begin{equation}
\begin{array}{llrll}
\dot x^2_5&-&x^4_5&=&\phantom{-}F_4(t)\\
\noalign{\ms}
\dot x^4_5&+&V_{22}(t)x^2_5&=&-F_2(t).\\
\end{array}
\end{equation}
The initial conditions are $x^r_5(0)=0$, as we see again from
(16). In particular we have
\begin{equation}
x^1_5(t)=x^3_5(t)=0,
\end{equation}
and since we have the fundamental matrix (57) of the
homogeneous system (54), we obtain the
the remaining two functions $x^2_5(t)$ and $x^4_5(t)$ by
{\it variation of constants} as
\begin{equation}
\left(
\begin{array}{l}
x^2_5(t)\\
\noalign{\sms}
x^4_5(t)\\
\end{array}
\right) =
\left(
\begin{array}{ll}
x^2_2(t)&x^2_4(t)\\
\noalign{\sms}
x^4_2(t)&x^4_4(t)\\
\end{array}
\right)
\int\limits_0^td\tau\left[
\left(
\begin{array}{rr}
x^4_4(\tau)&-x^2_4(\tau)\\
\noalign{\sms}
-x^4_2(\tau)&x^2_2(\tau)\\
\end{array}
\right)
\left(
\begin{array}{r}
F_4(\tau)\\
\noalign{\sms}
-F_2(\tau)\\
\end{array}
\right)
\right].
\end{equation}
Note that with (50) and (56)-(62) all flow derivatives $x^i_\lambda(t)$
of first order are now determined. Among them the $x^i_2(t)$ and
$x^i_0(t)$ have been seen to vanish for $i=1,3$, as announced at
the end of section~6, and so there is in fact no need to
determine third derivatives
$Z^2_{\ell m n}$ and $Z^0_{\ell m n}$ of the starting point and
flow time functions.

\vskip1.5cm

\noi
{\sefont 9 Calculation of higher order flow derivatives}

\bs\noi
Just as the first equation of (49) led to (53) and (54), so the
other equations of (49) show that the higher order flow
derivatives $x^r\st(t)$ satisfy differential equations
\begin{equation}
\begin{array}{llrll}
\dot x^1\st&-&x^3\st&=&g^1\st(t)\\
\noalign{\ms}
\dot x^3\st&+&V_{11}(t)x^1\st&=&g^3\st(t)\\
\end{array}
\end{equation}
and
\begin{equation}
\begin{array}{llrll}
\dot x^2\st&-&x^4\st&=&g^2\st(t)\\
\noalign{\ms}
\dot x^4\st&+&V_{22}(t)x^2\st&=&g^4\st(t),\\
\end{array}
\end{equation}
with
\begin{equation}
\begin{array}{lll}
g^r_{\lambda\mu}(t)&=&
v^r_{ij}(t)x^i_\lambda(t) x^j_\mu(t)\quad\text{and}\\
\noalign{\bs}
g^r_{\lambda\mu\nu}(t)&=&
v^r_{ijk}(t)x^i_\lambda(t) x^j_\mu(t) x^k_\nu(t)\\
\noalign{\ms}
&&+v^r_{ij}(t)\big(x^i_{\lambda}(t)x^j_{\mu\nu}(t)+
x^i_{\mu}(t)x^j_{\nu\lambda}(t)
+x^i_{\nu}(t)x^j_{\lambda\mu}(t)\big)\\
\end{array}
\end{equation}

\ms\noi
What are the initial conditions? Again from (16) we see that
\begin{equation}
x^i_{\lambda\mu}(0)=x^i_{\lambda\mu\nu}(0)=0\quad\text{if}\quad
\lambda,\mu,\nu\ne0
\end{equation}
and therefore by variation of constants we get
\begin{equation}
\left(
\begin{array}{l}
x^1\st(t)\\
\noalign{\sms}
x^3\st(t)\\
\end{array}
\right) =
\left(
\begin{array}{ll}
x^1_1(t)&x^1_3(t)\\
\noalign{\sms}
x^3_1(t)&x^3_3(t)\\
\end{array}
\right)
\int\limits_0^td\tau\left[
\left(
\begin{array}{rr}
x^3_3(\tau)&-x^1_3(\tau)\\
\noalign{\sms}
-x^3_1(\tau)&x^1_1(\tau)\\
\end{array}
\right)
\left(
\begin{array}{r}
g^1\st(\tau)\\
\noalign{\sms}
g^3\st(\tau)\\
\end{array}
\right)
\right]
\end{equation}
and
\begin{equation}
\left(
\begin{array}{l}
x^2\st(t)\\
\noalign{\sms}
x^4\st(t)\\
\end{array}
\right) =
\left(
\begin{array}{ll}
x^2_2(t)&x^2_4(t)\\
\noalign{\sms}
x^4_2(t)&x^4_4(t)\\
\end{array}
\right)
\int\limits_0^td\tau\left[
\left(
\begin{array}{rr}
x^4_4(\tau)&-x^2_4(\tau)\\
\noalign{\sms}
-x^4_2(\tau)&x^2_2(\tau)\\
\end{array}
\right)
\left(
\begin{array}{r}
g^2\st(\tau)\\
\noalign{\sms}
g^4\st(\tau)\\
\end{array}
\right)
\right]
\end{equation}
for all indices $*=\lambda\mu$ and $*=\lambda\mu\nu$ with
$\lambda,\mu,\nu\ne0$. But do we know the functions $g^i\st(t)$
in all these cases? Let us look at (65). The $v^r\st(t)$ are all
known, see (48) and and (51). The $x^i_\lambda(t)$ have been
determined in section 8, so we know
all $g^r_{\lambda\mu}(t)$
and hence the $x^r_{\lambda\mu}(t)$ for $\lambda,\mu\ne0$ from
(67) and (68). These in turn give us, now for
$\lambda,\mu,\nu\ne0$,
the $g^r_{\lambda\mu\nu}(t)$ by (65) and the  $x^r_{\lambda\mu\nu}(t)$
from (67) and (68).

It remains to determine the $x^r_{\lambda\mu}(t)$ and
$x^r_{\lambda\mu\nu}(t)$ in those cases where one or several of
the indices are zero. The values of these functions at
$t=T(\ve_0)$ might also be needed in (25) for the calculation of
the partial derivatives of the Poincar\'e map. The index 0 denotes
the time derivative. Knowing the $x^r_0(t)$ from (59) we derive
$x^r_{00}(t)$ and $x^r_{000}(t)$ as
\begin{equation}
\begin{array}{lcl}
x^1_{00}(t)&=&\phantom{-}0\\
\noalign{\sms}
x^2_{00}(t)&=&\phantom{-}\ddot y(t)=-V_2(t)\quad\text{by}\quad(8)\\
\noalign{\sms}
x^3_{00}(t)&=&\phantom{-}0\\
\noalign{\sms}
x^4_{00}(t)&=&-V_{22}(t)\dot y(t),\\
\noalign{\bs}
x^1_{000}(t)&=&\phantom{-}0\\
\noalign{\sms}
x^2_{000}(t)&=&-V_{22}(t)\dot y(t)\\
\noalign{\sms}
x^3_{000}(t)&=&\phantom{-}0\\
\noalign{\sms}
x^4_{000}(t)&=&-V_{222}(t)\dot y(t)^2+V_2(t)V_{22}(t).\\
\end{array}
\end{equation}
For $\lambda,\mu\ne0$ the $x^r_{0\lambda}(t)$ are obtained by
(53), (54) and (58), (60), (61) from the known $x^i_\lambda(t)$ and
similarly the $x^r_{0\lambda\mu}(t)$ by (63), (64) and the first
equation of (65) from the $x^i_{\lambda}(t)$ and
$x^i_{\lambda\mu}(t)$. Finally, differentiating (53), (54),
(58), (60) and (61) we see that for $\lambda\in\{1,2,3,4\}$
\begin{equation}
\begin{array}{lcl}
x^1_{00\lambda}(t)&=&-V_{11}(t)x^1_\lambda(t)\\
\noalign{\sms}
x^2_{00\lambda}(t)&=&-V_{22}(t)x^2_\lambda(t)\\
\noalign{\sms}
x^3_{00\lambda}(t)&=&-V_{112}(t)
\dot y(t)x^1_\lambda(t)+V_{11}(t)x^3_\lambda(t)\\
\noalign{\sms}
x^4_{00\lambda}(t)&=&-V_{222}(t)
\dot y(t)x^2_\lambda(t)+V_{22}(t)x^4_\lambda(t),\\
\end{array}
\end{equation}
in particular $x^r_{00\lambda}(t)=0$ for $r+\lambda$ odd and
$r,\lambda\in\{1,2,3,4\}$,
and
\begin{equation}
\begin{array}{lcl}
x^1_{005}(t)&=&\phantom{-}0\\
\noalign{\sms}
x^2_{005}(t)&=&-V_{22}(t)x^2_5(t)
+F_{42}(t)\dot y(t)+F_{44}(t)V_2(t)\\
\noalign{\sms}
x^3_{005}(t)&=&\phantom{-}0\\
\noalign{\sms}
x^4_{005}(t)&=&-V_{222}(t)
\dot y(t)x^2_5(t)+V_{22}(t)x^4_5(t)
-F_{22}(t)\dot y(t)-F_{24}(t)V_2(t).\\
\end{array}
\end{equation}
In principle we now have all we need to calculate the partial
derivatives of the Poincar\'e map.

\vskip1.5cm

\noi
{\sefont 10 Summary of the procedure}

\bs\noi
Once the numerical prerequisites of section~2 are established,
we get the $x^i_\lambda(t)$ as described in section~8 almost
without further calculation, the only exceptions are $x^2_5$ and
$x^4_5$, for which the integral (62) has to be evaluated. As
explained in section~9,
we also have the $x^i_{\lambda\mu}(t)$ and
$x^i_{\lambda\mu\nu}(t)$ in those cases where at most one of the
indices $\lambda,\mu,\nu$ is different from zero. Next determine
the $x^r_{\lambda\mu}(t)$ for $\lambda,\mu\ne0$ by calculating
$g^r_{\lambda\mu}(t)$ from the first equation of (65) and
applying (67) and (68). Then we get
$x^i_{0\lambda\mu}(t)$ for $\lambda,\mu\ne0$ by (63),(64) and
by the first equation of (65) without new integration. Also the
flow time derivatives $Z^0_{\ell m}$ of (42) are now known.
Finally, we now have the $g^r_{\lambda\mu\nu}(t)$ for
$\lambda,\mu,\nu\ne0$ from the second equation of (65) and
we can calculate the corresponding $x^r_{\lambda\mu\nu}(t)$ as
the integrals (67) and (68). Taking values at $t=T(\ve_0)$ of
all these functions and applying (25), we obtain the 38 partial
derivatives of the Poincar\'e map at $(0,0,\ve_0,0)$ we wanted.

For the computer, these instructions may be good enough, but a
person might want to see step by step what is going on.
For this we have some choice in which order to proceed.
We will first describe all those steps that are not
connected with the deformation question.

\vskip1.5cm

\noi
{\sefont 11 The undeformed system}

\bs\noi{\bf Step 1.} Choose the potential $V(x,y)$ to be
studied, with $\frac{\de V}{\de x}(0,y)\equiv0$, choose one of
its libration families on the $y$-axis and a reference point $E_0$ for the
energy parameter $\ve=E-E_0$. Choose a fixed $\ve_0$ at which
the bifurcation behavior of the libration shall be predicted.

\bs\noi{\bf Step 2.} Set up a first part of
the numerical prerequisites, namely
$y(t)$, $\xi_1(t)$, $\xi_2(t)$
and their first derivatives, as
described in section~2, including the period $T(\ve_0)$. Define
the four functions $x^r_\lambda(t)$ on $[0,T(\ve_0)]$ with
$r,\lambda\in\{1,3\}$ by (56).

\bs\noi{\bf Step 3.} Collect the Jacobian matrix of the Poincar\'e
map at $(0,0,\ve_0)$, or monodromy matrix of our librating
orbit, as
\begin{equation}
\left(
\begin{array}{ll}
Q_q&Q_p\\
\noalign{\sms}
P_q&P_p\\
\end{array}
\right)
=
\left(
\begin{array}{ll}
x^1_1(T(\ve_0))&x^1_3(T(\ve_0))\\
\noalign{\sms}
x^3_1(T(\ve_0))&x^3_3(T(\ve_0))\\
\end{array}
\right),
\end{equation}
according to (25), (26) and (55). If the trace $Q_q+P_p$ is
different from $+2$,
the fixed point is regular and the orbit will not bifurcate. In
this case the procedure may stop here, since then we might not
be interested in the higher derivatives. But $Q_q+P_p=2$ is not
a {\it technical} necessity for going on.

\bs\noi{\bf Step 4.} Now we will calculate
the $\ve$-derivative $\Tr'_A(\ve_0)=Q_{q\ve}+P_{p\ve}$
of the trace. From (25) we find
\begin{equation}
\begin{array}{lcr}
Q_{q\ve}&=&
\displaystyle
\phantom{-}\frac{1}{V_2V_2}x^3_1x^4_2+\frac{1}{V_2}x^1_{12}\\
\noalign{\bs}
P_{p\ve}&=&
\displaystyle-\frac{V_{11}}{V_2V_2}x^1_3x^4_2+\frac{1}{V_2}x^3_{32}\\
\end{array}
\end{equation}
at $(0,0,\ve_0)$. Here we need the remaining functions $\eta_1(t)$
and $\eta_2(t)$ of the prerequisites: they define the
$x^r_\lambda(t)$ on $[0,T(\ve_0)]$ with
$r,\lambda\in\{2,4\}$ by (57). Apart from the factor
$x^4_2(T(\ve_0))$ in the first summand, they are needed
as functions on $[0,T(\ve_0)]$ to
calculate $x^1_{12}$ and $x^3_{32}$ in the second summand by
integration (67), because the inhomogeneities $g^1_{\lambda\mu}(t)$
and $g^3_{\lambda\mu}(t)$ for $\lambda\in\{1,3\}$ and
$\mu\in\{2,4\}$ turn out by (65) to be
\begin{equation}
\begin{array}{lcl}
g^1_{\lambda\mu}(t)&=&\phantom{-}0\\
\noalign{\ms}
g^3_{\lambda\mu}(t)&=&-V_{112}(t)x^1_\lambda(t) x^2_\mu(t).\\
\end{array}
\end{equation}
If step~3 has shown $\ve_0$ to be singular ($\Tr_A(\ve_0)=2$),
then after completion of step~4 we know if it is a
{\it cross-bifurcation}, that is if $\Tr_A'(\ve_0)\ne0$.

\bs\noi{\bf Step 5.} Is this cross-bifurcation
transcritical? To answer this question, we need the monodromy
matrix (72) from step~3 and the second partial derivatives of
$Q$ and $P$ by the variable $q$ and $p$ at $(0,0,\ve_0)$,
that is the $(x^i\circ Z)_{\ell m}$ for $i\in\{1,3\}$ and
$\ell,m\in\{1,2\}$,
to see if $\widetilde P_{\widetilde q\widetilde q}\ne0$, where
the `tilde' denotes adapted coordinates.
By (26) and from our knowledge of the $x^i_\lambda$, the second
equation of (25) reads
\begin{equation}
\begin{array}{lcl}
Q_{qq}&=&x^1_{11}\\
\noalign{\ms}
Q_{qp}&=&x^1_{13}\\
\noalign{\ms}
Q_{pp}&=&x^1_{33}\\
\noalign{\ms}
P_{qq}&=&x^3_{11}\\
\noalign{\ms}
P_{qp}&=&x^3_{13}\\
\noalign{\ms}
P_{pp}&=&x^3_{33}\\
\noalign{\ms}
\end{array}
\end{equation}
at $t=T(\ve_0)$.
The $x^i_{\lambda\mu}(t)$ for $i,\lambda,\mu\in\{1,3\}$
have to be calculated from (67) by integration with inhomogeneities
\begin{equation}
\begin{array}{lcl}
g^1_{\lambda\mu}(t)&=&0\\
\noalign{\ms}
g^3_{\lambda\mu}(t)&=&-V_{111}(t)x^1_\lambda(t)x^1_\mu(t).\\
\end{array}
\end{equation}

\bs\noi{\bf Step 6.} If the cross-bifurcation is not
transcritical, then we are interested in
\begin{equation}
\eps_B''(0)=
\frac{3\widetilde Q_{\widetilde q\widetilde q}
\widetilde P_{\widetilde q\widetilde p}
-\widetilde Q_{\widetilde p}
\widetilde P_{\widetilde q\widetilde q\widetilde q}}
{3\widetilde Q_{\widetilde p}\widetilde P_{\widetilde q\eps}}\,,
\end{equation}
since the bifurcation is fork-like if and only if
$\eps_B''(0)\ne0$, and its sign and absolute value describe
geometric properties of the fork. We need information beyond the
first five steps only for
$\widetilde P_{\widetilde q\eps}$ and
$\widetilde P_{\widetilde q\widetilde q\widetilde q}$. In the
present step~6 we will take care of
$\widetilde P_{\widetilde q\eps}$. For this, we only have to
complete step~4 by the calculation of $Q_{p\ve}$ and $P_{q\ve}$,
which turn out to be
\begin{equation}
\begin{array}{lcr}
Q_{p\ve}&=&
\displaystyle
\phantom{-}\frac{1}{V_2V_2}x^3_3x^4_2+\frac{1}{V_2}x^1_{32}\\
\noalign{\bs}
P_{q\ve}&=&\displaystyle
-\frac{V_{11}}{V_2V_2}x^1_1x^4_2+\frac{1}{V_2}x^3_{12}\\
\end{array}
\end{equation}
similar to (73), with $x^1_{32}$ and $x^3_{12}$ determined by
integration (67) with inhomogeneities given by (74).

\bs\noi{\bf Step 7.} To calculate
$\widetilde P_{\widetilde q\widetilde q\widetilde q}$, we will
need the third order partial derivatives of $Q$ and $P$ by $q$
and $p$. As a preparation we will now extend steps~4 and~6 to
the determination of {\it all}\, $x^r_{\lambda\mu}(t)$ for
$r,\lambda,\mu\in\{1,2,3,4\}$ by (65) and (67), (68). In all
these cases we have
\begin{equation}
g^1_{\lambda\mu}(t)=g^2_{\lambda\mu}(t)=0
\end{equation}
and
\begin{equation}
\begin{array}{cll}
g^3_{\lambda\mu}(t)&=&\left\{
\begin{array}{cl}
-V_{111}(t)x^1_\lambda(t) x^1_\mu(t)&
\quad\text{if}\quad \lambda,\mu\in\{1,3\}\\
\noalign{\ms}
-V_{112}(t)x^1_\lambda(t) x^2_\mu(t)&
\quad\text{if}\quad \lambda\in\{1,3\}\;\text{and}\;\mu\in\{2,4\}\\
\noalign{\ms}
0&\quad\text{if}\quad \lambda,\mu\in\{2,4\}\\
\end{array}\right.\\
\noalign{\vskip1cm}
g^4_{\lambda\mu}(t)&=&\left\{
\begin{array}{cl}
-V_{112}(t)x^1_\lambda(t) x^1_\mu(t)&
\quad\text{if}\quad \lambda,\mu\in\{1,3\}\\
\noalign{\ms}
0&
\quad\text{if}\quad \lambda\in\{1,3\}\;\text{and}\;\mu\in\{2,4\}\\
\noalign{\ms}
-V_{222}(t)x^2_\lambda(t) x^2_\mu(t)&
\quad\text{if}\quad \lambda,\mu\in\{2,4\}\\
\end{array}\right.\\
\end{array}
\end{equation}
Note that therefore
\begin{equation}
\begin{array}{cl}
x^1_{\lambda\mu}(t)=x^3_{\lambda\mu}(t)=0&\quad\text{if}\quad
\lambda,\mu\in\{2,4\}\quad\text{and}\\
\noalign{\ms}
x^2_{\lambda\mu}(t)=x^4_{\lambda\mu}(t)=0&\quad\text{if}\quad
\lambda\in\{1,3\}\quad\text{and}\quad \mu\in\{2,4\}.\\
\end{array}
\end{equation}
To
calculate the other functions
$x^r_{\lambda\mu}(t)$ for
$r,\lambda,\mu\in\{1,2,3,4\}$ by (67) and (68) would be step~7.

\bs\noi{\bf Step 8.} It is now a suitable moment to improve upon
(42) and write out the second order flow time derivatives as
\begin{equation}
\begin{array}{lll}
Z^0_{11}&=&\displaystyle\frac1{V_2} x^4_{11}-\frac{V_{11}}{V_2V_2}x^4_2\\
\noalign{\ms}
Z^0_{12}&=&\displaystyle\frac1{V_2} x^4_{13}\\
\noalign{\ms}
Z^0_{22}&=&\displaystyle\frac1{V_2} x^4_{33}-\frac{1}{V_2V_2}x^4_2\\
\noalign{\ms}
Z^0_{13}&=&\;0\\
\noalign{\ms}
Z^0_{23}&=&\;0\,,\\
\noalign{\ms}
\end{array}
\end{equation}
again at $(0,0,\ve_0,0)$, with the $x^4\st$ to be taken
at $t=T(\ve_0)$ in the notation introduced with (47). We will
need these numbers when we apply the third equation of (25) to
calculate $Q_{qqq},Q_{qqp},\dots,P_{ppp}$. But before we can do
this, we have to determine the $x^r_{\lambda\mu\nu}(T(\ve_0))$
for $r,\lambda,\mu,\nu\in\{1,2,3,4\}$.

\break
\bs\noi{\bf Step 9.} This step is parallel to and based on step~7. We
read the inhomogeneities $g^r_{\lambda\mu\nu}(t)$ from the
second equation of (65) and then refer to (67), (68) to
determine the $x^r_{\lambda\mu\nu}(t)$ by integration. Again
\begin{equation}
g^1_{\lambda\mu\nu}(t)=g^2_{\lambda\mu\nu}(t)=0
\end{equation}
for all $\lambda,\mu,\nu\in\{1,2,3,4\}$, and the
$g^3_{\lambda\mu\nu}(t)$ are given in the four cases a)-d) as
follows.

\bs\noi
a) If $\lambda,\mu,\nu\in\{1,3\}$, then
\begin{equation}
\begin{array}{cll}
g^3_{\lambda\mu\nu}(t)&=&
-V_{1111}(t)x^1_\lambda(t) x^1_\mu(t)x^1_{\nu}(t)\\
\noalign{\ms}
&&-V_{111}(t)\big(
x^1_\lambda(t)x^1_{\mu\nu}(t)+
x^1_\mu(t)x^1_{\nu\lambda}(t)+
x^1_\nu(t)x^1_{\lambda\mu}(t)
\big)\\
\noalign{\ms}
&&-V_{112}(t)\big(
x^1_\lambda(t)x^2_{\mu\nu}(t)+
x^1_\mu(t)x^2_{\nu\lambda}(t)+
x^1_\nu(t)x^2_{\lambda\mu}(t)
\big)\,,\\
\end{array}
\end{equation}

\bs\noi
b) if $\lambda,\mu\in\{1,3\}$ and $\nu\in\{2,4\}$, then
\begin{equation}
\begin{array}{cll}
g^3_{\lambda\mu\nu}(t)&=&
-V_{1112}(t)x^1_\lambda(t) x^1_\mu(t)x^2_{\nu}(t)\\
\noalign{\ms}
&&-V_{111}(t)\big(
x^1_\lambda(t)x^1_{\mu\nu}(t)+
x^1_\mu(t)x^1_{\nu\lambda}(t)
\big)\phantom{\mkern8mu+
x^1_\nu(t)x^2_{\lambda\mu}(t)}\\
\noalign{\ms}
&&-V_{112}(t)
x^2_\nu(t)x^1_{\lambda\mu}(t)\,,
\\
\end{array}
\end{equation}

\bs\noi
c) if $\lambda\in\{1,3\}$ and $\mu,\nu\in\{2,4\}$, then
\begin{equation}
\begin{array}{cll}
g^3_{\lambda\mu\nu}(t)&=&
-V_{1122}(t)x^1_\lambda(t) x^2_\mu(t)x^2_{\nu}(t)\\
\noalign{\ms}
&&-V_{112}(t)\big(
x^1_\lambda(t)x^2_{\mu\nu}(t)+
x^2_\mu(t)x^1_{\nu\lambda}(t)+
x^2_\nu(t)x^1_{\lambda\mu}(t)
\big)\,,\\
\end{array}
\end{equation}

\bs\noi
d) if $\lambda,\mu,\nu\in\{2,4\}$, then
\begin{equation}
\begin{array}{cll}
g^3_{\lambda\mu\nu}(t)&=&\phantom{-}0.
\mkern350mu\\
\end{array}
\end{equation}

\bs\noi
Similarly the inhomogeneities $g^4_{\lambda\mu\nu}(t)$:

\ms\noi
a) If $\lambda,\mu,\nu\in\{1,3\}$, then
\begin{equation}
\begin{array}{cll}
g^4_{\lambda\mu\nu}(t)&=&
-V_{1112}(t)x^1_\lambda(t) x^1_\mu(t)x^1_{\nu}(t)\\
\noalign{\ms}
&&-V_{112}(t)\big(
x^1_\lambda(t)x^1_{\mu\nu}(t)+
x^1_\mu(t)x^1_{\nu\lambda}(t)+
x^1_\nu(t)x^1_{\lambda\mu}(t)
\big)\,,\\
\end{array}
\end{equation}

\bs\noi
b) if $\lambda,\mu\in\{1,3\}$ and $\nu\in\{2,4\}$, then
\begin{equation}
\begin{array}{cll}
g^4_{\lambda\mu\nu}(t)&=&
-V_{1122}(t)x^1_\lambda(t) x^1_\mu(t)x^2_{\nu}(t)\\
\noalign{\ms}
&&-V_{112}(t)\big(
x^1_\lambda(t)x^1_{\mu\nu}(t)+
x^1_\mu(t)x^1_{\nu\lambda}(t)
\big)\phantom{\mkern8mu+
x^1_\nu(t)x^2_{\lambda\mu}(t)}\\
\noalign{\ms}
&&-V_{222}(t)
x^2_\nu(t)x^2_{\lambda\mu}(t)\,,
\\
\end{array}
\end{equation}

\bs\noi
c) if $\lambda\in\{1,3\}$ and $\mu,\nu\in\{2,4\}$, then
\begin{equation}
\begin{array}{cll}
g^4_{\lambda\mu\nu}(t)&=&\phantom{-}0\,,
\mkern350mu\\
\end{array}
\end{equation}

\bs\noi
d) if $\lambda,\mu,\nu\in\{2,4\}$, then
\begin{equation}
\begin{array}{cll}
g^4_{\lambda\mu\nu}(t)&=&
-V_{2222}(t)x^2_\lambda(t) x^2_\mu(t)x^2_{\nu}(t)\\
\noalign{\ms}
&&-V_{222}(t)\big(
x^2_\lambda(t)x^2_{\mu\nu}(t)+
x^2_\mu(t)x^2_{\nu\lambda}(t)+
x^2_\nu(t)x^2_{\lambda\mu}(t)
\big).\\
\end{array}
\end{equation}
Note that in particular
\begin{equation}
\begin{array}{cl}
x^1_{\lambda\mu\nu}(t)=x^3_{\lambda\mu\nu}(t)=0&\quad\text{if}
\;\lambda,\mu,\nu\in\{2,4\}\\
\noalign{\ms}
x^2_{\lambda\mu\nu}(t)=x^4_{\lambda\mu\nu}(t)=0&\quad\text{if}
\;\lambda\in\{1,3\}\;\text{and}\; \mu,\nu\in\{2,4\}.
\end{array}
\end{equation}
The other functions
$x^r_{\lambda\mu\nu}(t)$ for
$r,\lambda,\mu,\nu\in\{1,2,3,4\}$
have to be calculated by integrations (67) and (68) in step~9.

\bs\noi{\bf Step 10.} We can now write down the eight missing
numbers $Q_{qqq},Q_{qqp},\dots,P_{ppp}$, thereby completing the
first part of our program, the part that is not concerned with deformation.
The $Z^2_{\ell m}$ and $Z^0_{\ell m}$ with $\ell,m\in\{1,2\}$ in
the formulas are taken from (36) and (82). The third
equation of (25) gives:
\begin{equation}
\begin{array}{lcrrrrrrl}
Q_{qqq}&=&x^1_{111}&+&3x^1_{12}Z^2_{11}&+&3x^3_1Z^0_{11}&&\\
\noalign{\ms}
Q_{qqp}&=&x^1_{113}&+&x^1_{32}Z^2_{11}&+&x^3_3Z^0_{11}
&+&2x^3_1Z^0_{12}\\
\noalign{\ms}
Q_{qpp}&=&x^1_{133}&+&x^1_{12}Z^2_{22}&+&x^3_1Z^0_{22}
&+&2x^3_1Z^0_{12}\\
\noalign{\ms}
Q_{ppp}&=&x^1_{333}&+&3x^1_{32}Z^2_{22}&+&3x^3_3Z^0_{22}&&\\
\noalign{\bs}
P_{qqq}&=&x^3_{111}&+&3x^3_{12}Z^2_{11}&-&3V_{11}x^1_1Z^0_{11}&&\\
\noalign{\ms}
P_{qqp}&=&x^3_{113}&+&x^3_{32}Z^2_{11}&-&V_{11}x^1_3Z^0_{11}
&-&2V_{11}x^1_1Z^0_{12}\\
\noalign{\ms}
P_{qpp}&=&x^3_{133}&+&x^3_{12}Z^2_{22}&-&V_{11}x^1_3Z^0_{22}
&-&2V_{11}x^1_1Z^0_{12}\\
\noalign{\ms}
P_{ppp}&=&x^3_{333}&+&3x^3_{32}Z^2_{22}&-&3V_{11}x^1_3Z^0_{22}&&\\
\noalign{\ms}
\end{array}
\end{equation}
at $t=T(\ve_0)$.

\vskip1.5cm

\noi
{\sefont 12 Deformation}

\bs\noi{\bf Step 11.} Choose a deformation term $F(x,y,p_x,p_y)$
satisfying the libration preserving condition (6). No new
`numerical prerequisites' are required, we can start right away
calculating $x^2_5(t)$ and $x^4_5(t)$ from (62), which is the
same as (68) with $g^2_5(t)=F_4(t)$ and $g^4_5(t)=-F_2(t)$. For
the other $x^i_5(t)$ see (50) and (61).

\bs\noi{\bf Step 12.} Now we can determine the sixteen functions
$x^r_{\lambda5}(t)$ for $r,\lambda\in\{1,2,3,4\}$. Again we
derive the corresponding $g^r_{\lambda5}(t)$ from (65). It turns
out that $g^r_{\lambda5}(t)=0$ if $r+\lambda$ is odd and hence
also
\begin{equation}
x^r_{\lambda5}(t)=0\quad\text{if}\quad
r+\lambda \quad\text{is odd.}
\end{equation}
For $\lambda\in\{1,3\}$ we get
\begin{equation}
\begin{array}{lcl}
g^1_{\lambda5}(t)&=&
\phantom{-}F_{13}(t)x^1_\lambda(t)+F_{33}(t)x^3_\lambda(t)\\
\noalign{\ms}
g^3_{\lambda5}(t)&=&
-F_{11}(t)x^1_\lambda(t)-F_{13}(t)x^3_\lambda(t)
-V_{112}(t)x^1_\lambda(t) x^2_5(t),\\
\end{array}
\end{equation}
and if
$\lambda\in\{2,4\}$, then
\begin{equation}
\begin{array}{lcl}
g^2_{\lambda5}(t)&=&
\phantom{-}F_{24}(t)x^2_\lambda(t)+F_{44}(t)x^4_\lambda(t)\\
\noalign{\ms}
g^4_{\lambda5}(t)&=&
-F_{22}(t)x^2_\lambda(t)-F_{24}(t)x^4_\lambda(t)
-V_{222}(t)x^2_\lambda(t) x^2_5(t).\\
\end{array}
\end{equation}

\bs\noi{\bf Step 13.} Now we collect the needed flow time and starting
point derivatives at $(0,0,\ve_0,0)$ that involve the
deformation parameter $u^4=\delta$, namely the $Z^\lambda_4$ and
$Z^\lambda_{\ell4}$ for $\lambda=0,2$ and $\ell=1,2$. From (36),
(41) and (42) we get
\begin{equation}
\begin{array}{lll}
Z^2_{4}&=&-\displaystyle\frac{F}{V_2}\\
\noalign{\ms}
Z^0_4\mkern4mu&=&\phantom{-}\displaystyle\frac1{V_2}
x^4_5-\frac{F}{V_2V_2}x^4_2\\
\noalign{\ms}
\end{array}
\end{equation}
and
\begin{equation}
\begin{array}{lll}
Z^2_{14}&=&\phantom{-}Z^2_{24}\quad =\quad 0\\
\noalign{\ms}
Z^0_{14}&=&\phantom{-}\displaystyle\frac{1}{V_2} x^4_{15}
\mkern110mu\\
\noalign{\ms}
Z^0_{24}&=&\phantom{-}\displaystyle\frac{1}{V_2} x^4_{35}\\
\noalign{\ms}
\end{array}
\end{equation}
at $(0,0,\ve_0,0)$.

\bs\noi{\bf Step 14.} We can now calculate those second
derivatives of the Poincar\'e map in which the deformation
parameter is involved:
\begin{equation}
\begin{array}{lcrrrrr}
Q_{q\delta}&=&x^1_{15}&+&x^1_{12}Z^2_{4}&+&x^3_1Z^0_4\\
\noalign{\ms}
Q_{p\delta}&=&x^1_{35}&+&x^1_{32}Z^2_{4}&+&x^3_3Z^0_{4}\\
\noalign{\ms}
P_{q\delta}&=&x^3_{15}&+&x^3_{12}Z^2_{4}&-&V_{11}x^1_1Z^0_{4}\\
\noalign{\ms}
P_{p\delta}&=&x^3_{35}&+&x^3_{32}Z^2_{4}&-&V_{11}x^1_3Z^0_{4}\\
\noalign{\ms}
\end{array}
\end{equation}
at $t=T(\ve_0)$.
But before we reach the {\it third}\, derivatives of $Q$ and $P$
involving $\delta$, we have to take one more step.


\bs\noi{\bf Step 15.}
We have to determine the
$x^r_{\lambda\mu5}(T(\ve_0))$ for $\lambda,\mu\in\{1,2,3,4\}$,
again by integration (67) and (68), with inhomogeneities
$g^r_{\lambda\mu5}(t)$ as follows.

\bs\noi
1a) If $\lambda,\mu,\in\{1,3\}$, then
\begin{equation}
\begin{array}{cll}
g^1_{\lambda\mu5}(t)&=&
\phantom{+}F_{113}(t)x^1_\lambda(t) x^1_\mu(t)\\
\noalign{\ms}
&&+F_{133}(t)\Big[x^1_\lambda(t) x^3_\mu(t)
+x^3_\lambda(t) x^1_\mu(t)\Big]\\
\noalign{\ms}
&&+F_{333}(t)x^3_\lambda(t) x^3_\mu(t)\\
\noalign{\ms}
&&+F_{13}(t)x^1_{\lambda\mu}(t)+F_{33}(t)x^3_{\lambda\mu}(t),\\
\end{array}
\end{equation}

\bs\noi
1b) if $\lambda\in\{1,3\}$ and $\mu\in\{2,4\}$, then
\begin{equation}
\begin{array}{cll}
g^1_{\lambda\mu5}(t)&=&
\phantom{+}F_{233}(t)x^3_\lambda(t) x^3_\mu(t)
+F_{334}(t)x^3_\lambda(t) x^4_\mu(t)\\
\noalign{\ms}
&&+F_{123}(t)x^1_\lambda(t) x^2_\mu(t)
+F_{134}(t)x^1_\lambda(t) x^4_\mu(t)\\
\noalign{\ms}
&&+F_{13}(t)x^1_{\lambda\mu}(t)+F_{33}(t)x^3_{\lambda\mu}(t),\\
\end{array}
\end{equation}

\bs\noi
1c) if $\lambda,\mu\in\{2,4\}$, then
\begin{equation}
g^1_{\lambda\mu5}(t)=0.
\end{equation}

\vskip1cm\noi
2a) If $\lambda,\mu,\in\{1,3\}$, then
\begin{equation}
\begin{array}{cll}
g^2_{\lambda\mu5}(t)&=&
\phantom{+}F_{114}(t)x^1_\lambda(t) x^1_\mu(t)\\
\noalign{\ms}
&&+F_{134}(t)\Big[x^1_\lambda(t) x^3_\mu(t)
+x^3_\lambda(t) x^1_\mu(t)\Big]\\
\noalign{\ms}
&&+F_{334}(t)x^3_\lambda(t) x^3_\mu(t)\\
\noalign{\ms}
&&+F_{24}(t)x^2_{\lambda\mu}(t)+F_{44}(t)x^4_{\lambda\mu}(t),\\
\end{array}
\end{equation}

\bs\noi
2b) if $\lambda\in\{1,3\}$ and $\mu\in\{2,4\}$, then
\begin{equation}
g^2_{\lambda\mu5}(t)=0,
\end{equation}

\bs\noi
2c) if $\lambda,\mu,\in\{2,4\}$, then
\begin{equation}
\begin{array}{cll}
g^2_{\lambda\mu5}(t)&=&
\phantom{+}F_{224}(t)x^2_\lambda(t) x^2_\mu(t)\\
\noalign{\ms}
&&+F_{244}(t)\Big[x^2_\lambda(t) x^4_\mu(t)
+x^4_\lambda(t) x^2_\mu(t)\Big]\\
\noalign{\ms}
&&+F_{444}(t)x^4_\lambda(t) x^4_\mu(t)\\
\noalign{\ms}
&&+F_{24}(t)x^2_{\lambda\mu}(t)+F_{44}(t)x^4_{\lambda\mu}(t).\\
\end{array}
\end{equation}

\vskip1cm\noi
3a) If $\lambda,\mu,\in\{1,3\}$, then
\begin{equation}
\begin{array}{cll}
g^3_{\lambda\mu5}(t)&=&
-V_{1112}(t)x^1_\lambda(t) x^2_\mu(t)x^2_{5}(t)\\
\noalign{\ms}
&&-V_{111}(t)\Big[
x^1_\lambda(t)x^1_{\mu5}(t)+
x^1_\mu(t)x^1_{\lambda5}(t)
\Big]\\
\noalign{\ms}
&&-V_{112}(t)\Big[
x^1_\lambda(t)x^2_{\mu5}(t)+
x^1_\mu(t)x^2_{\lambda5}(t)+
x^2_5(t)x^1_{\lambda\mu}(t)
\Big]\\
\noalign{\ms}
&&-F_{111}(t)x^1_\lambda(t) x^1_\mu(t)\\
\noalign{\ms}
&&-F_{113}(t)\Big[x^1_\lambda(t) x^3_\mu(t)
+x^3_\lambda(t) x^1_\mu(t)\Big]\\
\noalign{\ms}
&&-F_{133}(t)x^3_\lambda(t) x^3_\mu(t)\\
\noalign{\ms}
&&-F_{11}(t)x^1_{\lambda\mu}(t)-F_{13}(t)x^3_{\lambda\mu}(t),\\
\end{array}
\end{equation}

\bs\noi
3b) if $\lambda\in\{1,3\}$ and $\mu\in\{2,4\}$, then
\begin{equation}
\begin{array}{cll}
g^3_{\lambda\mu5}(t)&=&
-V_{1122}(t)x^1_\lambda(t) x^2_\mu(t)x^2_{5}(t)\\
\noalign{\ms}
&&-V_{122}(t)\Big[
x^1_\lambda(t)x^2_{\mu5}(t)+
x^2_\mu(t)x^1_{\lambda5}(t)+
x^2_5(t)x^1_{\lambda\mu}(t)
\Big]\\
\noalign{\ms}
&&-F_{112}(t)x^1_\lambda(t) x^2_\mu(t)
-F_{114}(t)x^1_\lambda(t) x^4_\mu(t)\\
\noalign{\ms}
&&+F_{123}(t)x^3_\lambda(t) x^2_\mu(t)
+F_{134}(t)x^3_\lambda(t) x^4_\mu(t)\\
\noalign{\ms}
&&-F_{11}(t)x^1_{\lambda\mu}(t)-F_{13}(t)x^3_{\lambda\mu}(t),\\
\end{array}
\end{equation}

\bs\noi
3c) if $\lambda,\mu\in\{2,4\}$, then
\begin{equation}
g^3_{\lambda\mu5}(t)=0.
\end{equation}

\vskip1cm\noi
4a) If $\lambda,\mu\in\{1,3\}$ then
\begin{equation}
\begin{array}{cll}
g^4_{\lambda\mu5}(t)&=&
-V_{1122}(t)x^1_\lambda(t) x^1_\mu(t)x^2_{5}(t)\\
\noalign{\ms}
&&-V_{112}(t)\Big[
x^1_\lambda(t)x^1_{\mu5}(t)+
x^1_\mu(t)x^1_{\lambda5}(t)
\Big]\\
\noalign{\ms}
&&-V_{222}(t)
x^2_5(t)x^2_{\lambda\mu}(t)\\
\noalign{\ms}
&&-F_{112}(t)x^1_\lambda(t) x^1_\mu(t)
-F_{233}(t)x^3_\lambda(t) x^3_\mu(t)\\
\noalign{\ms}
&&-F_{123}(t)\Big[x^1_\lambda(t) x^3_\mu(t)
+x^3_\lambda(t) x^1_\mu(t)\Big]\\
\noalign{\ms}
&&-F_{22}(t)x^2_{\lambda\mu}(t)-F_{24}(t)x^4_{\lambda\mu}(t),\\
\end{array}
\end{equation}

\bs\noi
4b) if $\lambda\in\{1,3\}$ and $\mu\in\{2,4\}$, then
\begin{equation}
g^4_{\lambda\mu5}(t)=-V_{222}(t)x^2_\mu(t)x^2_{\lambda5}(t),
\end{equation}

\bs\noi
4c) if $\lambda,\mu\in\{2,4\}$, then
\begin{equation}
\begin{array}{cll}
g^4_{\lambda\mu5}(t)&=&
-V_{2222}(t)x^2_\lambda(t) x^2_\mu(t)x^2_{5}(t)\\
\noalign{\ms}
&&-V_{222}(t)\Big[
x^2_\lambda(t)x^2_{\mu5}(t)+
x^2_\mu(t)x^2_{\lambda5}(t)+
x^2_5(t)x^2_{\lambda\mu}(t)
\Big]\\
\noalign{\ms}
&&-F_{222}(t)x^2_\lambda(t) x^2_\mu(t)\\
\noalign{\ms}
&&-F_{224}(t)\Big[x^2_\lambda(t) x^4_\mu(t)
+x^4_\lambda(t) x^2_\mu(t)\Big]\\
\noalign{\ms}
&&-F_{244}(t)x^4_\lambda(t) x^4_\mu(t)\\
\noalign{\ms}
&&-F_{22}(t)x^2_{\lambda\mu}(t)-F_{24}(t)x^4_{\lambda\mu}(t).\\
\end{array}
\end{equation}

\vskip1cm\noi
Note that in particular
\begin{equation}
x^1_{\lambda\mu5}(t)=x^3_{\lambda\mu5}(t)=0\quad\text{if}\quad
\lambda,\mu\in\{2,4\}.
\end{equation}

\bs\noi{\bf Step 16.} This then will be the last step in the
procedure to determine the 38 numbers we set out to calculate.
Only six of them remain, and these are given by
\begin{equation}
\begin{array}{lcrrrrrrl}
Q_{qq\delta}&=&x^1_{115}&+&x^1_{112}Z^2_{4}
&+&x^3_{11}Z^0_{4}&+&2x^3_1Z^0_{14}\\
\noalign{\ms}
Q_{qp\delta}&=&x^1_{135}&+&x^1_{123}Z^2_{4}
&+&x^3_{13}Z^0_{4}&+&x^3_1Z^0_{14}+x^3_3Z^0_{24}\\
\noalign{\ms}
Q_{pp\delta}&=&x^1_{335}&+&x^1_{233}Z^2_{4}
&+&x^3_{33}Z^0_{4}&+&2x^3_3Z^0_{24}\\
\noalign{\ms}
P_{qq\delta}&=&x^3_{115}&+&x^3_{112}Z^2_{4}&-&V_{11}x^1_{11}Z^0_{4}
&-&2V_{11}x^1_1Z^0_{14}\\
\noalign{\ms}
P_{qp\delta}&=&x^3_{135}&+&x^3_{123}Z^2_{4}&-&V_{11}x^1_{13}Z^0_{4}
&-&V_{11}(x^1_3Z^0_{14}+x^1_1Z^0_{24})\\
\noalign{\ms}
P_{pp\delta}&=&x^3_{335}&+&x^3_{233}Z^2_{4}&-&V_{11}x^1_{13}Z^0_{4}
&-&2V_{11}x^3_1Z^0_{24}\\
\noalign{\ms}
\end{array}
\end{equation}
at $t=T(\ve_0)$, as before, or explicitly: the left hand sides
are meant to be taken at $(q,p,\ve,\delta)=(0,0,\ve_0,0)$, thus
on the right hand sides the $x^r\st(t)$ at $t=T(\ve_0)$,
the $V\st(x,y)$
at $(x,y)=(0,y_{\rm max}(\ve_0))$, and the $Z^r\st$
are also to be taken at $(q,p,\ve,\delta)=(0,0,\ve_0,0)$,
as in
(41), (42), (82), (97), (98) where they have been calculated
from $x^r\st(T(\ve_0))$, $V\st(0,y_{\rm max}(\ve_0))$ and
various $F\st(0,y_{\rm max}(\ve_0),0,0)$.

\end{document}